\documentclass[]{interact}

\usepackage{epstopdf}
\usepackage[caption=false]{subfig}

\usepackage[numbers,sort&compress]{natbib}
\bibpunct[, ]{[}{]}{,}{n}{,}{,}
\makeatletter
\def\NAT@def@citea{\def\@citea{\NAT@separator}}
\makeatother
\usepackage{anyfontsize}
\theoremstyle{plain}
\newtheorem{theorem}{Theorem}[section]
\newtheorem{lemma}[theorem]{Lemma}
\newtheorem{corollary}[theorem]{Corollary}
\newtheorem{proposition}[theorem]{Proposition}

\theoremstyle{definition}
\newtheorem{definition}[theorem]{Definition}
\newtheorem{remark}[theorem]{Remark}
\newtheorem{example}[theorem]{Example}

\usepackage{xcolor}

\theoremstyle{remark}

\begin{document}
	

	\title{On the generalized Cauchy dual of closed operators in Hilbert spaces}
 \author{
		\name{Arup Majumdar\textsuperscript{a} \thanks{Arup Majumdar (corresponding author). Email address: arupmajumdar93@gmail.com},  P. Sam Johnson\textsuperscript{b}, Ram N. Mohapatra\textsuperscript{c}}
		\affil{\textsuperscript{a,}\textsuperscript{b}Department of Mathematical and Computational Sciences,
			National Institute of Technology Karnataka (NITK), Surathkal, Mangaluru 575025, India.} \textsuperscript{c} Department of Mathematics, University of Central Florida, Orlando, FL. 32816, USA}
  
	\maketitle

	\begin{abstract}
		In this paper, we introduce the generalized Cauchy dual $w(T) = T(T^{*}T)^{\dagger}$ of a closed operator $T$ with the closed range between Hilbert spaces and present intriguing findings that characterize the Cauchy dual of $T$. Additionally, we establish the result $w(T^{n}) = (w(T))^{n}$, for all $n \in \mathbb{N}$, where $T$ is a quasinormal EP operator.
	\end{abstract}
	
      \begin{amscode}47A05, 47B02, 47B20.\end{amscode}
    \section{Introduction}
The concept of the Cauchy dual of left-invertible operators was introduced by Shimorin \cite{MR1810120} as part of his work on Wold-type decompositions and wandering subspaces. Specifically, the Cauchy dual of a bounded right-invertible operator $T^{*}$ (the adjoint of $T$) is defined as the operator $T^{\prime} = T(T^{*}T)^{-1}$. In 2007, Chavan further explored the Cauchy dual in the context of $2$-hyperexpansive operators \cite{MR2360521}. Since then, there have been additional developments in the study of the Cauchy dual for both bounded and unbounded operators \cite{MR3906434, MR2784021}. In 2019, H. Ezzahraoui et al. examined interesting properties of the Cauchy dual for bounded closed-range operators in Hilbert spaces \cite{MR3967887}. In that paper, the authors defined the Cauchy dual of a bounded closed-range operator $T$ as $w(T) = T(T^{*}T)^{\dagger}$, where $(T^{*}T)^{\dagger}$ is the Moore-Penrose inverse of the operator $T^{*}T$. In this work, we introduce the term ``generalized Cauchy dual" to refer to this operator $w(T)$, as we employ the generalized inverse of $T^{*}T$ instead of the usual inverse $(T^{*}T)^{-1}$. We further explore the characterizations of the generalized Cauchy dual for closed operators in Hilbert spaces in Section 2.
    
From now on, we shall mean $H$, $K$, $H_{i}$, $K_{i}$ ($i = 1,2,\dots, n$) as Hilbert spaces. The specification of a domain is an essential part of the definition of an unbounded operator, usually defined on a subspace.  Consequently, for an operator $T$, the specification of the subspace $D$ on which $T$ is defined, called the domain of $T$, denoted by $D(T)$, is to be given. The null space and range space of $T$ are denoted by $N(T)$ and $R(T)$, respectively. $W_{1}^{\perp}$ denotes the orthogonal complement of a set $W_{1}$ whereas $W_{1} \oplus^{\perp} W_{2}$ denotes the orthogonal direct sum of the subspaces $W_{1}$ and $W_{2}$ of a Hilbert space. Moreover,  $T_{\vert_{W}}$ denotes the restriction of $T$ to a subspace $W$ of a specified Hilbert space.  We call $D(T)\cap N(T)^\perp $, the carrier of $T$ and it is denoted by $ C(T)$. $T^{*}$ denotes the adjoint of $T$, when $D(T)$ is densely defined in the specified Hilbert space. Here, $P_{V}$ is the orthogonal projection on the closed subspace $V$ in the specified Hilbert space and the set of bounded operators from $H$ into $K$ is denoted by $B(H, K)$. Similarly, the set of all bounded operators on $H$ is denoted by $B(H)$. Furthermore, $I_{W}$ stands for the identity operator on a subspace $W$ of the specified Hilbert space. 
	For the sake of completeness of exposition, we first begin with the definition of a closed operator.

	\begin{definition}
		Let $ T$ be an operator from a Hilbert space $H$ with domain $D(T) $ to a Hilbert space $K$. If the graph of $T$ defined by 
		$$ G(T)=\left\{(x,Tx): x\in D(T)\right\} $$ is closed in $H\times K $, then $T$ is called a closed  operator. Equivalently, $T $ is a closed operator if $ \{x_n \}$ in $D(T) $ such that $ x_n\rightarrow x $ and $ Tx_n\rightarrow y$ for some $ x\in  H,y\in   K $,  then $ x\in  D(T) $ and $ Tx=y $. That is, $G(T)$ is  a closed subspace of $H\times K$ with respect to the graph norm $\|(x, y)\|_T=(\|x\|^2+\|y\|^2)^{1/2}$. It is easy to show that the graph norm  $\|(x, y)\|_T$ is equivalent to the norm $\|x\|+\|y\|$.  We note that,  for any densely defined closed operator $T$,  the closure of $C(T)$, that is,  $ \overline{C(T)}$ is $ N(T)^\perp$. $C(H)$ denotes the set of all closed operators from $H$ into $H$. Meanwhile, $C(H, K)$ denotes the set of all closed operators from $H$ into $K$.
		
    \end{definition}

 We say that $S$ is an extension of $ T $ (denoted by $ T\subset S$) if $ D(T)\subset D(S) $ and $ Sx=Tx $ for all $ x\in D(T)$. Moreover, $T \in C(H)$ commutes $A \in B(H)$  if $AT \subset TA$.
	\begin{definition}
		Let $T$ be a closed operator from  $D(T) \subset H$ to  $K$. The Moore-Penrose inverse of $T$ is the map $T^{\dagger}: R(T) \oplus^{\perp} R(T)^{\perp} \to H$ defined by
		\begin{equation}\label{equ 1}
			T^{\dagger} y = 
			\begin{cases}
				({T}\vert_{C(T)})^{-1}y  & \text{if} ~  y\in R(T)\\
				0    & \text{if}  ~ y\in R(T)^{\perp}.
			\end{cases}
		\end{equation}
    \end{definition}
\noindent  It can be shown that $T^{\dagger}$ is bounded if and only if $R(T)$ is closed, when $T$ is closed.
	
 \begin{definition}\cite{MRARUP}\label{def 1.6}
 Let $T$ be a densely defined closed operator on $H$. The operator $T$ is said to be an EP operator if $T$ has closed range and $R(T) =R(T^{*})$.\\
 Moreover, the densely defined closed operator $T \in C(H)$ is called a hypo-EP operator if $R(T)$ is closed with $R(T) \subset R(T^{*})$.
 \end{definition}
\begin{theorem}\cite{MRARUP}\label{thm 2.6.4}
Let $T$ be a closed EP operator. Then $T^{n}$ is also EP and $R(T^{n}) = R(T)$, for all $n \in \mathbb{N}$.
\end{theorem}
 \begin{theorem}\label{thm 1.7}\cite{MR0396607}
  Let $T$ be a densely defined closed operator from $D(T) \subset H$ into $K$. Then the following statements hold:
  \begin{enumerate}
  \item $T^{\dagger}$ is a closed operator from $K$ into $H$;
  \item $D(T^{\dagger}) = R(T) \oplus^{\perp} N(T^{*})$; $N(T^{\dagger}) = N(T^{*})$;
  \item $R(T^{\dagger}) = C(T)$;
  \item $T^{\dagger}Tx = P_{\overline{R(T^{\dagger})}}x, \text{ for all } x\in D(T)$;
  \item $TT^{\dagger}y = P_{\overline{R(T)}}y, \text{ for all } y\in D(T^{\dagger})$;
  \item $(T^{\dagger})^{\dagger} = T$;
  \item $(T^{*})^{\dagger} = (T^{\dagger})^{*}$;
  \item $N((T^{*})^{\dagger})= N(T)$;
  \item $(T^{*}T)^{\dagger} = T^{\dagger}(T^{*})^{\dagger}$;
  \item $(TT^{*})^{\dagger} = (T^{*})^{\dagger} T^{\dagger}$.
  \end{enumerate}
  \end{theorem}
   The space $H \bigoplus K$ defined by $H \bigoplus K = \{(h, k): h\in H, k\in K\}$ is a linear space with respect to addition and scalar multiplication defined by 
\begin{align*}
&~~~~~~~~~~~~~~~~~~~~~~(h_{1}, k_{1}) + (h_{2}, k_{2}) = (h_{1} + h_{2}, k_{1} + k_{2}), \text{ and }\\
&\lambda (h, k ) = (\lambda h, \lambda k), \text{ for all } h, h_{1}, h_{2}\in H, \text{ for all } k, k_{1}, k_{2}\in K \text{ and } \lambda\in \mathbb{K}, ~(\mathbb{K} = \mathbb{R} \text{ or } \mathbb{C}).
\end{align*}
Now, $H \bigoplus K$ is an inner product space with respect to the inner product given by
\begin{align*}
\langle (h_{1}, k_{1}) (h_{2}, k_{2}) \rangle = \langle h_{1}, h_{2}\rangle + \langle k_{1}, k_{2}\rangle, \text{ for all } h_{1}, h_{2}\in H, \text{ and } \text{ for all }k_{1}, k_{2} \in K. 
\end{align*}
The norm on $H \bigoplus K$ is defined by 
\begin{align*}
\|(h, k)\| = (\|h\|^{2} + \|k\|^{2})^{\frac{1}{2}}, \text{ for all } (h,k) \in H \times K.
\end{align*}
Moreover, the direct sum of two operators $T_{1} \text{ and } T_{2}$ from $D(T_{1}) \subset H_{1}$ to $K_{1}$ and from $D(T_{2}) \subset H_{2}$ to $K_{2}$ respectively is defined by 
\begin{align*}
(T_{1} \bigoplus T_{2}) (h_{1}, h_{2}) = (T_{1}h_{1}, T_{2}h_{2}), \text{ for all } h_{1} \in D(T_{1}), \text { and for all } h_{2}\in D(T_{2}).
\end{align*}
 
 \begin{theorem}\label{thm 1.8} \cite{MR007}
Let $T_{1}: D(T_{1}) \subset H_{1} \to K_{1}$ and $T_{2}: D(T_{2}) \subset  H_{2} \to K_{2}$ be two closed operators with closed ranges. Then $T = T_{1} \bigoplus T_{2} : D(T_{1}) \bigoplus D(T_{2}) \subset H_{1} \bigoplus H_{2} \to K_{1} \bigoplus K_{2}$ has the Moore-Penrose inverse. Moreover,
 \begin{align*}
 T^{\dagger} = (T_{1} \bigoplus T_{2})^{\dagger} = T_{1}^{\dagger} \bigoplus T_{2}^{\dagger}.
 \end{align*}
\end{theorem}
\begin{definition}
Let $T \in C(H, K)$ be a densely defined operator. The generalized Cauchy dual of $T$ is defined by $w(T) = T(T^{*}T)^{\dagger}$.
\end{definition} 
\begin{example}\label{example 1.8.8}
Define $T$ on $\ell^{2}$ by $$T(x_{1}, x_{2},\ldots,x_{n},\ldots)= (x_{2}, 0, 2x_{4},0, 3x_{6},\ldots,0,nx_{2n},0,\ldots)$$ with domain $D(T) = \{(x_{1}, x_{2},\ldots,x_{n},\ldots)\in \ell^{2} : \sum_{n=1}^\infty |n x_{2n}|^2<\infty	 \}$. Then $T$ is closed with closed range \cite{majumdar2024hyers}. Since, $N(T) = \overline{\text{span}\{e_{2n-1}: n \in \mathbb{N}\}}$, then
\begin{align*}
w(T)e_{2n-1}= (T^{*})^{\dagger}e_{2n-1}= 0  \text{ and } w(T)e_{2n}= (T^{*})^{\dagger}e_{2n}=\frac{1}{n}{e_{2n-1}}, \text{ for all } n \in \mathbb{N}.
\end{align*}
\end{example}
\begin{example}\label{example 1.9.8}
 Let $\mathcal{H}=L^{2}([a, b])$ and $K(.,.) \in L^{2}([a, b] \times[a, b])$.

$$
T f(x)=\int_{a}^{b} K(x, y) f(y) d y, \quad f \in H
$$

Suppose that $K(x, y)=\overline{K(y, x)}$. Then $T$ is compact and self-adjoint $\left(T=T^{*}\right)$. It is well known that in this case there exists $\left(\lambda_{n}\right)_{n \geq 1} \subset \mathbb{R}$ and an orthonormal sequence $\left(\varphi_{n}\right)_{n \geq 1} \subset L^{2}([a, b])$ such that $T=\sum_{n \geq 1} \lambda_{n} \varphi_{n} \otimes \varphi_{n}$. If $0 \neq \lambda_{k} \in \sigma(T)$, then $T-\lambda_{k} I$ is a bounded closed range operator and

$$
\omega\left(T-\lambda_{k} I\right)=\sum_{n \geq 1, n \neq k} \frac{1}{\lambda_{n}-\lambda_{k}} \varphi_{n} \otimes \varphi_{n}.
$$
\end{example}
\begin{definition}\cite{MR3192032}\label{def 1.7.8}
Let $T \in C(H)$ be a densely defined operator. $T$ is quasinormal if $T^{*}TT = TTT^{*}$.
\end{definition}
\begin{theorem}\cite{MR3192032}\label{thm 1.8.8}
Let $T \in C(H)$ be a densely defined quasinormal operator. Then $(T^{*}T)^{n} = (T^{*})^{n} T^{n} = (T^{n})^{*}T^{n}$, for all $n \in \mathbb{N}$.
\end{theorem}
\begin{theorem}\cite{MR3192032}\label{thm 1.9.8}
Let $T \in C(H)$ be a quasinormal operator. Then $T^{n}$ is also quasinormal, for all $n \in \mathbb{N}$.
\end{theorem}
 \section{Characterizations of the generalized Cauchy dual of closed operators in the Hilbert spaces}
 Throughout the section, we consider $T$ as the densely defined closed operator from  $H$ into $K$ with $R(T)$ to be closed. The properties of $w(T)$ will lead us to find the interesting characterizations of $(T^{*})^{\dagger}$ because $w(T) = T(T^{*}T)^{\dagger} = TT^{\dagger}(T^{*})^{\dagger} = (T^{*})^{\dagger}$.
 \begin{proposition}\label{pro 2.1.8}
 Let $T \in C(H, K)$ be a densely defined closed range operator. Then the following statements hold:
 \begin{enumerate}
 \item $w(T) = (T^{\dagger})^{*}$ is bounded.
 \item $w(w(T))$ is closable and $\overline{w(w(T))} = T$.
 \item $(w(T))^{*} = w(T^{*})$.
 \item $T^{*}w(T) = \overline{w(T^{*})T} = P_{R(T^{*})} $.
 \item $\overline{w(T)T^{*}} = Tw(T^{*}) = P_{R(T)}$.
 \item $w(T)^{*}w(T) = w(T^{*}T)$.
  \end{enumerate}
 \end{proposition}
 \begin{proof}
 $\textit{(1)}$ Since, $R(T)$ is closed which implies $T^{\dagger}$ is bounded, So, $w(T)= (T^{*})^{\dagger} = (T^{\dagger})^{*}$ is also bounded.\\
 $\textit{(2)}$ $w(w(T)) = w((T^{*})^{\dagger}) = (T^{*})^{\dagger} (T^{\dagger}(T^{*})^{\dagger})^{\dagger} = (T^{*})^{\dagger}(T^{*}T) = {P_{R(T)}}\vert_{D(T^{*})}T \subset T$. So, $w(w(T))$ is closable. Again $D(w(w(T)))$ is densely defined because $D(T^{*}T)$ is densely defined. Thus, $(w(w(T)))^{*} = (T^{*}T)T^{\dagger} = T^{*}P_{R(T)} = T^{*}$. Therefore, $\overline{w(w(T))} = T$.\\
 $\textit{(3)}$ By the closed range theorem, we get $R(T)$ is closed if and only if $R(T^{*})$ is closed. So, $w(T^{*}) = T^{\dagger}$. Thus, $(w(T))^{*}= T^{\dagger} = w(T^{*})$.\\
 $\textit{(4)}$  $T^{*}w(T) = (T^{*}T)(T^{*}T)^{\dagger} = P_{R(T^{*})}$. Moreover, $w(T^{*})T = T^{\dagger}T= P_{R(T^{*})}\vert_{D(T)} \subset T^{*}w(T)$. Hence, $\overline{w(T^{*})T} = T^{*}w(T) = P_{R(T^{*})}$.\\
 $\textit{(5)}$ $Tw(T^{*}) = TT^{\dagger} = P_{R(T)}$ and $w(T)T^{*} = P_{R(T)}\vert_{D(T^{*})}$ is closable with the minimal closed extension $P_{R(T)}$. Therefore, $\overline{w(T)T^{*}} = P_{R(T)} = Tw(T^{*})$.\\
 $\textit{(6)}$ $w(T)^{*}w(T) = T^{\dagger}(T^{\dagger})^{*} = (T^{*}T)^{\dagger} = w(T^{*}T)$.
 \end{proof}
 \begin{theorem}\label{thm 2.2.8}
 Let $T \in C(H)$ be a densely defined closed range operator. Then the following conditions are equivalent:
 \begin{enumerate}
 \item $T$ is selfadjoint.
 \item $T = \overline{w(T^{*})T}T^{*}$.
 \item $T^{*} = TT^{*}w(T)$.
 \end{enumerate}
 \end{theorem}
 \begin{proof}
 $\textit{(1) $\implies$ (2)}$ From the fourth statement of Proposition \ref{pro 2.1.8}, we have $\overline{w(T^{*})T}T^{*} = T^{*} = T$.\\
  $\textit{(2) $\implies$ (3)}$ Since, $\overline{w(T^{*})T}T^{*} = T$. From the third statement of Proposition \ref{pro 2.1.8}, we get $T^{*} = T (w(T^{*})T)^{*} = TT^{*} w(T)$.\\
  $\textit{(3) $\implies$ (1)}$ Now, $T^{*} = TT^{*}w(T) = TT^{*}(T^{*})^{\dagger} = TP_{R(T^{*})} = TP_{{N(T)}^{\perp}} = T$.
  \end{proof}
  \begin{theorem}\label{thm 2.3.8}
  Let $T \in C(H)$ be a densely defined closed range operator. Then $T$ is selfadjoint if and only if $w(T)$ is so. 
  \end{theorem}
  \begin{proof}
  Since $T = T^{*}$, then $(w(T))^{*} = w(T^{*}) = w(T)$. Hence, $w(T)$ is also selfadjoint.\\
  Conversely, 
 \begin{align*}
  (w(T))^{*} &= w(T)\\
  T^{\dagger} &= (T^{*})^{\dagger}\\
  (T^{\dagger})^{\dagger} &= ((T^{*})^{\dagger})^{\dagger}\\
  T &= T^{*}.
  \end{align*}
  Therefore, $T$ is selfadjoint.
  \end{proof}
\begin{theorem}\label{thm 2.4.8}
Let $T\in C(H)$ be a densely defined closed operator. Then $T$ is normal if and only if $w(T)$ is so.
\end{theorem}
\begin{proof}
Since, $TT^{*} = T^{*}T$. Then, from the third and sixth statements of Proposition \ref{pro 2.1.8}, we get $(w(T))^{*}w(T) = w(T^{*}T) = w(TT^{*})= w(T) (w(T))^{*}$. Hence, $w(T)$ is normal.\\
Conversely, 
\begin{align*}
w(T)(w(T^{*}))^{*} &= (w(T^{*}))^{*} w(T)\\
(T^{*})^{\dagger} T^{\dagger} &= T^{\dagger} (T^{*})^{\dagger}\\
((TT^{*})^{\dagger})^{\dagger} &= ((T^{*}T)^{\dagger})^{\dagger}\\
TT^{*} &= T^{*}T.
\end{align*}
Therefore, $T$ is normal.
\end{proof}
We now investigate the polar decomposition of $w(T)$. The polar decomposition of $T$ is represented as $T = U_{T}\vert T \vert$, where $\vert T \vert = (T^{*}T)^{\frac{1}{2}}$ and $U_{T}$ is partial isometry with the initial space $R(\vert T\vert)$ and the final space $R(T)$. 
\begin{lemma}\label{lem 2.5.8}
Let $T \in C(H, K)$ be a densely defined closed range operator. Then $T^{\dagger} = {\vert T \vert}^{\dagger} U_{T^{*}}$ and $(T^{*})^{\dagger} = U_{T}{\vert T \vert}^{\dagger}$.
\end{lemma}
\begin{proof}
We know the relations $R(U_{T}) = R(T), N(U_{T}) = N(T)= N(\vert T \vert)$, and $R(U_{T^{*}}) = R(T^{*}) = R(\vert T \vert),~ N(U_{T^{*}}) = N(T^{*})$. Now, we claim that $U_{T}U_{T^{*}} = P_{R(T)}$. Let us consider $x \in R(T)^{\perp} = N(T^{*})$, so $U_{T}U_{T^{*}}x = 0$. When $u \in R(T)$, then there exists $z \in C(T)$ such that $u = Tz$. Again,
$U_{T}U_{T^{*}}u = U_{T}\vert T\vert z = Tz = u$. Thus, $U_{T}U_{T^{*}}= P_{R(T)}$ and $U_{T^{*}}U_{T} = P_{R(T^{*})} = P_{R(\vert T\vert)}$. Hence, 
\begin{align}\label{equ 2.8}
T({\vert T \vert}^{\dagger} U_{T^{*}}) T = U_{T}\vert T \vert {\vert T \vert}^{\dagger} U_{T^{*}} U_{T} \vert T \vert = U_{T}\vert T \vert = T.
\end{align}
\begin{align}\label{equ 3.8}
({\vert T \vert}^{\dagger} U_{T^{*}})T({\vert T \vert}^{\dagger} U_{T^{*}})= {\vert T \vert}^{\dagger} U_{T^{*}}U_{T}\vert T \vert {\vert T \vert}^{\dagger} U_{T^{*}} = {\vert T \vert}^{\dagger} U_{T^{*}}.
\end{align}
\begin{align}\label{equ 4.8}
T({\vert T \vert}^{\dagger} U_{T^{*}}) = U_{T} \vert T\vert {\vert T \vert}^{\dagger} U_{T^{*}} = U_{T}U_{T^{*}} = P_{R(T)}.
\end{align}
and 
\begin{align}\label{equ 5.8}
({\vert T \vert}^{\dagger} U_{T^{*}})T = {\vert T \vert}^{\dagger} U_{T^{*}}U_{T}\vert T \vert =  {\vert T \vert}^{\dagger} \vert T \vert = P_{N(T)^{\perp}}\vert_{D(\vert T \vert)}.
\end{align}
By Theorem 5.7 \cite{MR0451661} and the relations (\ref{equ 2.8}), (\ref{equ 3.8}), (\ref{equ 4.8}), (\ref{equ 5.8}) confirm that $T^{\dagger} = {\vert T \vert }^{\dagger} U_{T^{*}}$.
Therefore, $(T^{*})^{\dagger} = U_{T}{\vert T\vert}^{\dagger}$.
\end{proof}
\begin{theorem}\label{thm 2.6.8}
Let $T \in C(H, K)$ be a densely defined closed range operator. Then $w(T) = U_{T}{\vert T \vert}^{\dagger}$ and $\|T^{\dagger}\| = \| {\vert T \vert}^{\dagger}\|$.
\end{theorem}
\begin{proof}
By Proposition \ref{pro 2.1.8}, $w(T) = (T^{*})^{\dagger}$ is bounded. So, $\|w(T)\| = \|T^{\dagger}\|$. Again,
$w(T) = TT^{\dagger}(T^{*})^{\dagger} = U_{T}\vert T \vert {\vert T \vert}^{\dagger} U_{T^{*}}U_{T} {\vert T \vert}^{\dagger} = U_{T} P_{R(T^{*})}{\vert T \vert}^{\dagger} = U_{T} {\vert T \vert}^{\dagger}$. Now, $\|w(T)x\| = \|U_{T}{\vert T \vert}^{\dagger} x\| = \|{\vert T \vert}^{\dagger}x\|$ because $R({\vert T \vert}^{\dagger}) \subset N(T)^{\perp} = R(T^{*}) = R(\vert T \vert)$. Hence, $\|w(T)\| = \|{\vert T \vert}^{\dagger}\| = \|T^{\dagger}\|$.
\end{proof}
\begin{remark}\label{remark 2.7.8}
When $T \in C(H)$ is a positive selfadjoint operator with closed range. Then $(T^{\dagger})^{\frac{1}{2}} = (T^{\frac{1}{2}})^{\dagger}$ (by using Theorem \ref{thm 1.7}$\mathit(9)$).
\end{remark}
\begin{remark}\label{remark 2.8.8}
When $T \in C(H)$ is a selfadjoint closed range operator then $w(T^{2}) = (w(T))^{2}$ (by using the sixth and third statements of Proposition \ref{pro 2.1.8}).
\end{remark}
\begin{corollary}\label{cor 2.9.8}
Let $T \in C(H, K)$ be a densely defined closed operator. Then $\vert{w(T)}\vert = w(\vert{T}\vert)$.
\end{corollary}
\begin{proof}
By Remark \ref{remark 2.8.8} and the third and sixth statements of Proposition \ref{pro 2.1.8}, we get $(w(\vert T \vert))^{2} = w(T^{*}T)= {\vert{w(T)}\vert}^{2}$. It is easy to show that $w(\vert T \vert) = {\vert T \vert}^{\dagger}$ is non-negative selfadjoint. Therefore, $\vert{w(T)}\vert = w(\vert{T}\vert)$.
\end{proof}
\begin{theorem}\label{thm 2.10.8}
Let $T \in C(H, K)$ be a densely defined closed range operator. Then $U_{T} = \overline{w(T)\vert T \vert}$ and the polar decomposition of $w(T) = U_{T} \vert{w(T)}\vert$.
\end{theorem}
\begin{proof}
From Theorem \ref{thm 2.6.8}, we have $w(T) = U_{T}{\vert{T}\vert}^{\dagger}$. So, $w(T)\vert T \vert = U_{T} (P_{{N(\vert T \vert)}^{\perp}}\vert_{D(\vert T \vert)}) \subset U_{T} P_{R(T^{*})}$. Thus, $w(T)\vert T \vert$ is closable and $(U_{T} P_{R(T^{*})})^{*} = (w(T)\vert T \vert)^{*}$. Hence, $\overline{w(T)\vert T \vert} = U_{T} P_{R(T^{*})} = U_{T}$.\\
Moreover, $w(T) = U_{T} {\vert T \vert}^{\dagger} = U_{T}w(\vert T \vert) = U_{T} \vert {w(T)}\vert$ (by Corollary \ref{cor 2.9.8}).
\end{proof}
\begin{lemma}\label{lemma 2.11.8}
Let $T \in C(H, K)$ be a densely defined closed range operator. Then $T\vert T \vert$ is a closed operator and $R(T \vert T \vert) = R(T)$.
\end{lemma}
\begin{proof}
We know $T\vert T\vert = U_{T}(T^{*}T)$. Let us consider an element $(x,y) \in \overline{G(T\vert T \vert)}$, then there exists a sequence $\{x_{n}\}$ in $D(T\vert T \vert)$ such that $x_{n} \to x$ and $U_{T}(T^{*}T)x_{n} \to y$ as $n \to \infty$. So, there is an element $z \in R(\vert T \vert)= R(T^{*})$ with $y = U_{T}z$. Again,
\begin{align*}
\|T^{*}Tx_{n}-z\| = \|U_{T}(T^{*}Tx_{n}-z)\| \to 0, \text{ as } n \to \infty.
\end{align*}
Since, $T^{*}T$ is closed which implies $x \in D(T^{*}T)$ and $T^{*}Tx = z$. Thus, $y = U_{T}T^{*}Tx = T\vert T \vert x$. Hence, $T\vert T \vert$ is closed. Moreover, $R(T\vert T \vert) = TR(T^{*}T) = TR(T^{*}) = R(TT^{*}) = R(T)$ is closed.
\end{proof}
\begin{remark}\label{remark 2.12.8}
When $T \in C(H, K)$ is a densely defined closed range operator, then $\vert T \vert T^{*} = (T^{*}T)U_{T^{*}}$ is closed with the closed range $R(\vert T \vert T^{*}) = T^{*}T R(U_{T^{*}}) = T^{*}TR(T^{*}) = R(T^{*}T) = R(T^{*})$.
\end{remark}
\begin{lemma}\label{lemma 2.13.8}
Let $T \in C(H, K)$ be a densely defined closed operator. Then $(T\vert T\vert)^{\dagger} = (\vert T \vert)^{\dagger} T^{\dagger}$.
\end{lemma}
\begin{proof}
By Lemma \ref{lemma 2.11.8}, we have $T \vert T \vert$ is closed with closed range, so $(T\vert T\vert)^{\dagger}$ is bounded. Now,
\begin{align}\label{equ 6.8}
\vert T \vert^{\dagger}T^{\dagger} T \vert T \vert \vert T \vert^{\dagger}T^{\dagger} = \vert T \vert^{\dagger} T^{\dagger}TT^{\dagger}= \vert T \vert^{\dagger}T^{\dagger}.
\end{align}
\begin{align}\label{equ 7.8}
T\vert T \vert  {\vert T \vert}^{\dagger}T^{\dagger}T {\vert T \vert} = TT^{\dagger}T\vert T \vert = T \vert T \vert.
\end{align}
\begin{align}\label{equ 8.8}
T\vert T \vert {\vert T \vert}^{\dagger}T^{\dagger} = TT^{\dagger} = P_{R(T)} \text{ is symmetric}.
\end{align}
\begin{align}\label{equ 9.8}
{\vert T \vert}^{\dagger}T^{\dagger}T \vert T \vert = ({\vert T \vert}^{\dagger}\vert T \vert)\vert_{D(T\vert T \vert)} = P_{R(T^{*})}\vert_{D(T^{*}T)} \text{ is symmetric because $D(T^{*}T)$ is dense}.
\end{align}
Therefore, Theorem 5.7 in \cite{MR0451661} and the relations (\ref{equ 6.8}), (\ref{equ 7.8}), (\ref{equ 8.8}), (\ref{equ 9.8}) confirm that  $(T\vert T\vert)^{\dagger} = (\vert T \vert)^{\dagger} T^{\dagger}$.
\end{proof}
\begin{lemma}\label{lemma 2.14.8}
Let $T \in C(H, K)$ be a densely defined closed range operator. Then $(\vert T\vert T^{*})^{\dagger} = (T^{*})^{\dagger} {\vert T \vert}^{\dagger}$.
\end{lemma}
\begin{proof}
By Remark \ref{remark 2.12.8}, we get $(\vert T \vert T^{*})^{\dagger}$ exists and is bounded. Now,
\begin{align}\label{equ 10.8}
(T^{*})^{\dagger}{\vert T \vert}^{\dagger} \vert T \vert T^{*}(T^{*})^{\dagger}\vert T \vert^{\dagger} = (T^{*})^{\dagger} \vert T \vert^{\dagger} \vert T \vert {\vert T \vert}^{\dagger} = (T^{*})^{\dagger}{\vert T \vert}^{\dagger}.
\end{align}
\begin{align}\label{equ 11.8}
\vert T \vert T^{*}(T^{*})^{\dagger} {\vert T \vert}^{\dagger} \vert T \vert T^{*} = \vert T \vert {\vert T \vert}^{\dagger} \vert T \vert T^{*}= \vert T \vert T^{*}.
\end{align}
\begin{align}\label{equ 12.8}
\vert T \vert T^{*}(T^{*})^{\dagger}{\vert T \vert}^{\dagger} = \vert T \vert {\vert T \vert}^{\dagger} = P_{R(T^{*})} \text{ is symmetric}.
\end{align}
\begin{align}\label{equ 13.8}
(T^{*})^{\dagger}{\vert T \vert}^{\dagger} \vert T \vert T^{*} = (T^{*})^{\dagger}T^{*}\vert_{D({\vert T \vert}T^{*})} = P_{R(T)}\vert_{D({\vert T \vert} T^{*})}.
\end{align}
Again, $D({\vert T \vert}T^{*}) \supset D({\vert T \vert}^{2}T^{*}) =D(T^{*}T T^{*}) \supset D(TT^{*}TT^{*})= D((TT^{*})^{2})$. We know $D((TT^{*})^{2})$ is dense in $D(TT^{*})$ and $D(TT^{*})$ is also dense in $K$. So, $D(\vert T \vert T^{*})$ is dense in $K$. Hence, the operator $(T^{*})^{\dagger}{\vert T \vert}^{\dagger} \vert T \vert T^{*} = P_{R(T)}\vert_{D({\vert T \vert} T^{*})}$ is symmetric. Therefore, Theorem 5.7 in \cite{MR0451661} and the relations (\ref{equ 10.8}), (\ref{equ 11.8}), (\ref{equ 12.8}), (\ref{equ 13.8}) guarantee that $(\vert T\vert T^{*})^{\dagger} = (T^{*})^{\dagger} {\vert T \vert}^{\dagger}$.
\end{proof}
\begin{theorem}\label{thm 2.15.8}
Let $T \in C(H)$ be a densely defined closed range operator. Then $\vert T\vert T^{*} = T \vert T \vert$ if and only if $\vert w(T)\vert w(T)^{*} = w(T) \vert w(T) \vert$.
\end{theorem}
\begin{proof}
\begin{align*}
\vert w(T) \vert w(T)^{*} &= {\vert T \vert}^{\dagger} T^{\dagger} ~(\text{by Corollary \ref{cor 2.9.8}})\\
&= (T\vert T \vert)^{\dagger} ~(\text{by Lemma \ref{lemma 2.13.8}})\\
&= (\vert T \vert T^{*})^{\dagger}\\
&= (T^{*})^{\dagger}{\vert T \vert}^{\dagger} ~(\text{by Lemma \ref{lemma 2.14.8}})\\
&= w(T)\vert w(T) \vert.
\end{align*}
Conversely,
\begin{align*}
\vert w(T) \vert w(T)^{*} &= w(T) \vert w(T) \vert\\
{\vert T \vert}^{\dagger} T^{\dagger} &= (T^{*})^{\dagger} \vert T \vert^{\dagger}\\
(T\vert T \vert)^{\dagger} &= (\vert T \vert T^{*})^{\dagger} ~(\text{by Lemma \ref{lemma 2.13.8}, and Lemma \ref{lemma 2.14.8}})\\
T\vert T \vert &= \vert T \vert T^{*}.
\end{align*}
\end{proof}
\begin{lemma}\label{lemma 2.16.8}
Let $T \in C(H)$ be a densely defined operator. If $T$ is quasinormal EP operator, then $T(T^{*}T)^{\dagger} \supset (T^{*}T)^{\dagger}T$.
\end{lemma}
\begin{proof}
Since, $T$ is EP then $R(T) = R(T^{*})$. So, we can represent $T$ as:
\[T =
\left[\begin{array}{cc}
T_{1} & 0 \\
0 & 0
\end{array}\right]:\left[\begin{array}{c}
C(T) \\
N(T)
\end{array}\right] \rightarrow\left[\begin{array}{c}
R(T) \\
N(T)
\end{array}\right], \text{ where $T_{1} = T\vert_{C(T)} : C(T) \to R(T)$}
.\]\\
Moreover, $N({T_{1}}^{*}) = \{0\}$. Again, ${T_{1}}^{-1}$ exists and is bounded. So, $({T_{1}}^{-1})^{*} = ({T_{1}}^{*})^{-1}$ is also bounded. Now, $D(T^{*}T) = D({T_{1}}^{*}T_{1}) \oplus^{\perp} N(T)$ and 
$T^{*}T =
\left[\begin{array}{cc}
{T_{1}}^{*}T_{1} & 0 \\
0 & 0
\end{array}\right]
.~ \text{Thus, $R({T_{1}}^{*}T_{1}) = R(T^{*}T)$ is closed.}$ Since, $T$ is quasinormal which implies $T_{1}(T_{1}^{*}T_{1}) = (T_{1}^{*}T_{1})T_{1}$. Furthermore, $({T_{1}}^{*}T_{1})^{\dagger} = {T_{1}}^{\dagger} ({T_{1}}^{*})^{\dagger} = {T_{1}}^{-1} (T_{1}^{*})^{-1}$. From the relation $T_{1}(T_{1}^{*}T_{1}) = (T_{1}^{*}T_{1})T_{1}$, we get $(T_{1}^{*})^{-1}T_{1}(T_{1}^{*}T_{1}) = I_{D(T_{1}^{*})} T_{1}^{2}$ which implies $(T_{1}^{*})^{-1}T_{1}T_{1}^{*} = I_{D(T_{1}^{*})} T_{1}$. Hence, $(T_{1}^{*})^{-1}T_{1} = I_{D(T_{1}^{*})}T_{1}(T_{1}^{*})^{-1}$. We know that $N(T_{1}^{*}T_{1}) = N(T_{1}) = \{0\}$. So, $(T_{1}^{*}T_{1})^{-1}$ exists and $(T_{1}^{*} T_{1})^{-1} = T_{1}^{-1} (T_{1}^{*})^{-1}$ because $(T_{1}^{*} T_{1}) T_{1}^{-1} (T_{1}^{*})^{-1}= I_{R(T)}$ and $T_{1}^{-1} (T_{1}^{*})^{-1} (T_{1}^{*}T_{1}) \subset I_{C(T)} \subset I_{N(T)^{\perp}}$. By Theorem \ref{thm 1.8}, we get
\[(T^{*}T)^{\dagger} =
\left[\begin{array}{cc}
(T_{1}^{*}T_{1})^{\dagger} & 0 \\
0 & 0
\end{array}\right]
.\]\\ So,
\[(T^{*}T)^{\dagger}T =
\left[\begin{array}{cc}
(T_{1}^{*}T_{1})^{\dagger}T_{1} & 0 \\
0 & 0
\end{array}\right]  =
\left[\begin{array}{cc}
T_{1}^{-1}(T_{1}^{*})^{-1}T_{1} & 0 \\
0 & 0
\end{array}\right]\]\\    \[(T^{*}T)^{\dagger}T = \left[\begin{array}{cc}
T_{1}^{-1}I_{D(T_{1}^{*})}T_{1}(T_{1}^{*})^{-1}) & 0 \\
0 & 0
\end{array}\right] \subset \left[\begin{array}{cc}
T_{1}^{-1}T_{1}(T_{1}^{*})^{-1}) & 0 \\
0 & 0
\end{array}\right]
.\]\\ Therefore, 
\[(T^{*}T)^{\dagger}T  \subset \left[\begin{array}{cc}
T_{1}^{-1}T_{1}(T_{1}^{*})^{-1} & 0 \\
0 & 0
\end{array}\right]  \subset \left[\begin{array}{cc}
T_{1}T_{1}^{-1}(T_{1}^{*})^{-1} & 0 \\
0 & 0
\end{array}\right] = \left[\begin{array}{cc}
T_{1}(T_{1}^{*}T_{1})^{\dagger} & 0 \\
0 & 0
\end{array}\right] = T(T^{*}T)^{\dagger}
.\]\\
\end{proof}
\begin{lemma}\label{lemma 2.17.8}
Let $T \in C(H, K)$ be a densely defined closed range operator. Then $(T^{*}T + P_{N(T)})$ is invertible and $(T^{*}T + P_{N(T)})^{-1} = T^{\dagger}(T^{*})^{\dagger} + P_{N(T)}$. Moreover, $T^{\dagger} = \overline{(T^{*}T + P_{N(T)})^{-1}T^{*}}$ and $w(T) = T(T^{*}T + P_{N(T)})^{-1}$.
\end{lemma}
\begin{proof}
It is obvious that $D(T^{*}T + P_{N(T)}) = D(T^{*}T)$. Let us consider $x \in N(T^{*}T + P_{N(T)})$, then $x = x_{1} + x_{2}$, where $x_{1} \in N(T^{*}T) = N(T),$ $x_{2} \in N(T)^{\perp} \cap D(T^{*}T)$ and $(T^{*}T + P_{N(T)})x = 0$ which implies $T^{*}Tx_{2} = -x_{1} \in N(T) \cap N(T)^{\perp} = \{0\}$. So, $x_{1} = 0$ and $x_{2} \in N(T) \cap N(T)^{\perp} = \{0\}$. Thus, $x=0$ which confirms the existence of $(T^{*}T + P_{N(T)})^{-1}$. Now, we claim that $(T^{*}T + P_{N(T)})^{-1} = T^{\dagger}(T^{*})^{\dagger} + P_{N(T)}$.
\begin{align}
(T^{*}T + P_{N(T)}) (T^{\dagger}(T^{*})^{\dagger} + P_{N(T)}) = T^{*}(T^{*})^{\dagger} + P_{N(T)} = P_{{N(T)}^{\perp}} + P_{N(T)}= I_{H}.
\end{align}
and 
\begin{align}
&(T^{\dagger}(T^{*})^{\dagger} + P_{N(T)}) (T^{*}T + P_{N(T)})\\ &= T^{\dagger}(T^{*})^{\dagger}T^{*}T + P_{N(T)}\\ &= {P_{{N(T)}^{\perp}}}\vert_{D(T^{*}T)} + P_{N(T)} \\ 
&= I_{D(T^{*}T)} \subset I_{H}.
\end{align}
Thus, $(T^{*}T + P_{N(T)})^{-1} = T^{\dagger}(T^{*})^{\dagger} + P_{N(T)}$ is bounded. This confirms the invertibility of $(T^{*}T + P_{N(T)})$. Again,
\begin{align*}
(T^{*}T + P_{N(T)})^{-1}T^{*} = T^{\dagger}(T^{*})^{\dagger}T^{*} = T^{\dagger}P_{R(T)}\vert_{D(T^{*})} \subset T^{\dagger}P_{R(T)} = T^{\dagger}.
\end{align*}
Hence, $(T^{\dagger})^{*} = T (T^{*}T + P_{N(T)})^{-1} = w(T)$. Therefore, $T^{\dagger} = \overline{(T^{*}T + P_{N(T)})^{-1}T^{*}}$.
\end{proof}
\begin{theorem}\label{thm 2.18.8}
Let $T \in C(H)$ be a quasinormal EP operator. Then $w(T^{n}) =(w(T))^{n}$, for all $n \in \mathbb{N}$.
\end{theorem}
\begin{proof}
Corollary 3.8 \cite{MR3192032} and Theorem \ref{thm 2.6.4} confirm that $T^{n}$ is quasinormal EP operator with $R(T^{n}) = R(T)$, for all $n \in \mathbb{N}$. We claim that $N(T^{n}) = N(T)$, for all $n \in \mathbb{N}$. It is obvious that $N(T) \subset N(T^{n})$, for all $n \in \mathbb{N}$. Let $x \in N(T^{2})$, then $Tx \in N(T) \cap N(T)^{\perp} = \{0\}$. So, $x \in N(T)$ which implies $N(T) = N(T^{2})$. Again, consider $z \in N(T^{3})$, then $T^{*}T^{3}z = 0$. By the definition of quasinormal, we get $TT^{*}T^{2}z = 0$ which implies $T^{2}z \in N(TT^{*}) \cap N(T)^{\perp} = N(T) \cap N(T)^{\perp} = \{0\}$. Thus, $z \in N(T^{2})$. It guarantees $N(T^{3}) = N(T^{2}) = N(T)$. Similarly, by the induction hypothesis, we have $N(T^{n}) = N(T)$, for all $n \in \mathbb{N}$. By Lemma \ref{lemma 2.17.8} and Theorem \ref{thm 1.8.8}, we get 
\begin{align*}
w(T^{n}) &= T^{n}((T^{n})^{*}T^{n} + P_{N(T^{n})})^{-1}\\
&= T^{n}((T^{n})^{\dagger} (({T^{n}})^{*})^{\dagger} + P_{N(T)})\\
&= T^{n}(({(T^{n})}^{*}T^{n})^{\dagger}  + P_{N(T)})\\
&= T^{n}(((T^{*}T)^{n})^{\dagger}).
\end{align*}
Moreover, $N(T^{*}T) = N(T) = N(T^{n}) = N((T^{n})^{*}T^{n}) = N((T^{*}T)^{n})$. Now we will show that $((T^{*}T)^{n})^{\dagger} = ((T^{*}T)^{\dagger})^{n}$, for all $n \in \mathbb{N}$. 
\begin{align*}
((T^{*}T)^{\dagger})^{n} (T^{*}T)^{n} &= (T^{*}T)^{\dagger}(T^{*}T)\vert_{D((T^{*}T)^{n})}\\ &= P_{{N(T^{*}T)}^{\perp}}\vert_{D((T^{*}T)^{n})} \\ &=P_{N((T^{*}T)^{n})^{\perp}}\vert_{D((T^{*}T)^{n})}.
\end{align*}
Furthermore,
\begin{align*}
(T^{*}T)^{n}((T^{*}T)^{\dagger})^{n} = (T^{*}T)(T^{*}T)^{\dagger} = P_{R(T^{*})} &= P_{R(T)} ~( \text{ since, $T$ is EP})\\
&= P_{R(T^{n})}~(\text{ By Theorem \ref{thm 2.6.4}})\\ 
&= P_{R((T^{n})^{*})} ~(\text{ since, $T^{n}$ is also EP})\\
&= P_{R((T^{n})^{*}(T^{n}))}\\
&= P_{R((T^{*}T)^{n})}.
\end{align*}
By Theorem 5.7 \cite{MR0451661}, we get $((T^{*}T)^{n})^{\dagger} = ((T^{*}T)^{\dagger})^{n}$.
\begin{align*}
(w(T))^{n} &= (T(T^{\dagger}(T^{*})^{\dagger} + P_{N(T)}))^{n}\\
&= (TT^{\dagger}(T^{*})^{\dagger})^{n}\\
&= (T(T^{*}T)^{\dagger})^{n}\\
& \subset T^{n}((T^{*}T)^{\dagger})^{n} (\text{ by Lemma \ref{lemma 2.16.8}})\\
&=  T^{n}((T^{*}T)^{n})^{\dagger} = w(T^{n}).
\end{align*}
We know that $(w(T))^{n}$ is bounded. Therefore, $w(T^{n}) =(w(T))^{n}$, for all $n \in \mathbb{N}$.
\end{proof}
\begin{remark}
In general, $w(T^{n}) = (w(T))^{n}$ is not true even in finite-dimensional Hilbert spaces. Consider $T = \left[\begin{array}{cc}
1 & 0 \\
1 & 0
\end{array}\right]$, then the generalized Cauchy dual $w(T) = \left[\begin{array}{cc}
\frac{1}{2} & 0 \\
\frac{1}{2} & 0
\end{array}\right] = \frac{1}{2}{T}$ and $T^{2} =T$. So, $w(T^{2}) = w(T) \neq (w(T))^{2} = \frac{1}{4}{T}$.
\end{remark}
\begin{remark}\label{remark 2.19.8}
 Theorem \ref{thm 2.18.8} is a generalized version of Remark \ref{remark 2.8.8} because every selfadjoint closed range operator is always quasinormal EP. So, Remark \ref{remark 2.8.8} can be generalized for all $n \in \mathbb{N}.$
\end{remark}
In general, $w(ST) \neq w(S)w(T)$, where $S \in B(H)$ and $T \in C(H)$ is a densely defined closed range operator. Because when $T \in C(H)$ is a densely defined closed range operator, then $T^{\dagger}$ is bounded. Then $w(T^{\dagger}) w(T) = ((T^{\dagger})^{*})^{\dagger} (T^{*})^{\dagger} = T^{*}(T^{*})^{\dagger} = P_{R(T^{*})}$. However, $T^{\dagger} T = P_{N(T)^{\perp}}\vert_{D(T)}$ is not closed but closable. Moreover, $w(\overline{T^{\dagger}T}) = w(P_{N(T)^{\perp}})= P_{N(T)^{\perp}} = P_{R(T^{*})} = w(T^{\dagger}) w(T)$. Additionally, $w(TT^{\dagger}) = w(P_{R(T)}) = P_{R(T)} = \overline{P_{N(T^{*})^{\perp}\vert_{D(T^{*})}}} = \overline{(T^{*})^{\dagger} T^{*}} = \overline{w(T)w(T^{\dagger})}$. In Theorem \ref{thm 2.20.8}, we will present a sufficient condition to obtain the relation $w(ST) = w(S)w(T)$.
\begin{theorem}\label{thm 2.20.8}
Let $T \in C(H)$ be a densely defined closed range operator and $S \in B(H)$ be an EP operator. If $R(S) = R(T)$, then $w(ST) = w(S)w(T)$.
\end{theorem}
\begin{proof}
Since, $R(S) = R(T)$. Then $T$ and $S$ can be written as:
\[T =
\left[\begin{array}{cc}
T_{1} & 0 \\
0 & 0
\end{array}\right]:\left[\begin{array}{c}
C(T) \\
N(T)
\end{array}\right] \rightarrow\left[\begin{array}{c}
R(S) \\
N(S)
\end{array}\right], \text{ where $T_{1} = T\vert_{C(T)} : C(T) \to R(S)$}
.\]\\
\[S =
\left[\begin{array}{cc}
S_{1} & 0 \\
0 & 0
\end{array}\right]:\left[\begin{array}{c}
R(S) \\
N(S)
\end{array}\right] \rightarrow\left[\begin{array}{c}
R(S) \\
N(S)
\end{array}\right], \text{ where $S_{1} : R(S) \to R(S)$}
.\]\\
Then $R(ST) = SR(T) = R(S^{2}) =R(S)$ is closed (by Theorem \ref{thm 2.6.4}). Again, $S_{1}T_{1}$ is closed because $T_{1}$ is closed and $S_{1}^{-1}$ is bounded. Thus, $ST$ is closed. Now, $w(ST) = ((ST)^{*})^{\dagger} = (T^{*}S^{*})^{\dagger}$ and $w(S)w(T) = (S^{*})^{\dagger} (T^{*})^{\dagger}$. We will show that $(T^{*}S^{*})^{\dagger} = (S^{*})^{\dagger} (T^{*})^{\dagger}$. 
\begin{align}\label{equ 19.8}
(S^{*})^{\dagger}(T^{*})^{\dagger}(T^{*}S^{*})(S^{*})^{\dagger}(T^{*})^{\dagger} = (S^{*})^{\dagger} (T^{*})^{\dagger}T^{*}P_{R(S)}(T^{*})^{\dagger} = (S^{*})^{\dagger} (T^{*})^{\dagger}.
\end{align}
\begin{align}\label{equ 20.8}
(T^{*}S^{*}) (S^{*})^{\dagger}(T^{*})^{\dagger}(T^{*}S^{*})= T^{*}(T^{*})^{\dagger}T^{*}S^{*} = T^{*}S^{*}.
\end{align}
\begin{align}\label{equ 21.8}
(T^{*}S^{*})(S^{*})^{\dagger}(T^{*})^{\dagger} = T^{*}(T^{*})^{\dagger} = P_{R(T^{*})} \text{ is symmetric}.
\end{align}
$D(T^{*}S^{*}) = D(ST)^{*}$ is dense because $ST$ is closed. Hence,
\begin{align}\label{equ 22.8}
(S^{*})^{\dagger}(T^{*})^{\dagger}T^{*}S^{*} = (S^{*})^{\dagger} P_{R(T)}\vert_{D(T^{*})}S^{*} = (S^{*})^{\dagger}S^{*}\vert_{D(T^{*}S^{*})} \subset P_{R(S)} \text{ is also symmetric}.
\end{align}
By Theorem 5.7 in \cite{MR0451661} and relations (\ref{equ 19.8}), (\ref{equ 20.8}), (\ref{equ 21.8}), (\ref{equ 22.8}), we have $(T^{*}S^{*})^{\dagger} = (S^{*})^{\dagger}(T^{*})^{\dagger}$. Therefore, $w(ST) = w(S)w(T)$.
\end{proof}
\begin{proposition}
Let $T_{i} \in C(H_{i}, K_{i}) ~(i = 1,2)$ be a densely defined closed range operator. Then $w(T_{1} \bigoplus T_{2}) = w(T_{1}) \bigoplus w(T_{2})$.
\end{proposition}
\begin{proof}
It is obvious that $T_{1} \bigoplus T_{2} $ is closed with closed range $R(T_{1} \bigoplus T_{2})$ (by Theorem \ref{thm 1.8}). Therefore,
\begin{align*}
w(T_{1} \bigoplus T_{2}) &= ((T_{1} \bigoplus T_{2})^{*})^{\dagger}\\ 
&= (T_{1}^{*} \bigoplus T_{2}^{*})^{\dagger}\\
&= (T_{1}^{*})^{\dagger} \bigoplus (T_{2}^{*})^{\dagger}  (\text{ by Theorem \ref{thm 1.8}})\\
&= w(T_{1}) \bigoplus w(T_{2}).
\end{align*}
\end{proof}
\begin{theorem}\label{thm 2.22.8}
\noindent Let $T =
\left[\begin{array}{cc}
T_{1} & T_{2} \\
0 & 0
\end{array}\right]$ be an operator on $H \bigoplus K$ with $T_{1} \in C(H)$ and $T_{2}\in C(K, H)$ both are densely defined closed range operators. If $R(T_{1}) \oplus^{\perp} R(T_{2}) \subset H$, then  $T^{\dagger} =
\left[\begin{array}{cc}
T_{1}^{\dagger} & 0 \\
T_{2}^{\dagger} & 0
\end{array}\right]$.
\end{theorem}
\begin{proof}
First, we will show that $T$ is a closed operator with closed range. Consider
$\biggl( \begin{pmatrix} x_{1}\\
x_{2}\end{pmatrix}, \begin{pmatrix} y_{1}\\y_{2} \end{pmatrix} \biggr )\in \overline{G(T)}$, then there exists a sequence $\biggl\{ \begin{pmatrix} x_{1n} \\ x_{2n} \end{pmatrix} \biggr\}$ such that $\biggl\{ \begin{pmatrix} x_{1n} \\ x_{2n} \end{pmatrix} \biggr\} \to \begin{pmatrix} x_{1}\\
x_{2}\end{pmatrix}$  and $\biggl\{ \begin{pmatrix} T_{1}x_{1n} + T_{2}x_{2n}\\ 0 \end{pmatrix} \biggr\} \to  \begin{pmatrix} y_{1}\\
y_{2}\end{pmatrix}$ as $n \to \infty$. We can observe that $y_{2} =0$ and $\{(T_{1}x_{1n} + T_{2}x_{2n})\}$ is a Cauchy sequence. The given condition $R(T_{1}) \oplus^{\perp} R(T_{2})$ says that $\{T_{1}x_{1n}\}$ and $\{T_{2}x_{2n}\}$ both are Cauchy sequences. Since, $T_{1}$ and $T_{2}$ both are closed. So, $T_{1}x_{1n} + T_{2}x_{2n} \to T_{1}x_{1} + T_{2}x_{2}$, as $n \to \infty$. Thus, $y_{1} = T_{1}x_{1} + T_{2}x_{2}$. Hence, $T$ is closed. $R(T) = R(T_{1}) \oplus^{\perp} R(T_{2})$ is closed because $R(T_{1}), R(T_{2})$ both are closed. So, $T^{\dagger}$ exists and is bounded. Now,
\begin{align}\label{equ 23.8}
 \left[\begin{array}{cc}
T_{1}^{\dagger} & 0 \\
T_{2}^{\dagger} & 0
\end{array}\right]
\left[\begin{array}{cc}
T_{1} & T_{2} \\
0 & 0
\end{array}\right] 
\left[\begin{array}{cc}
T_{1}^{\dagger} & 0 \\
T_{2}^{\dagger} & 0
\end{array}\right] =  \left[\begin{array}{cc}
T_{1}^{\dagger} & 0 \\
T_{2}^{\dagger} & 0
\end{array}\right]
\left[\begin{array}{cc}
T_{1}T_{1}^{\dagger} + T_{2}T_{2}^{\dagger} & 0 \\
0 & 0
\end{array}\right] = \left[\begin{array}{cc}
T_{1}^{\dagger} & 0 \\
T_{2}^{\dagger} & 0
\end{array}\right].
\end{align}
\begin{align}\label{equ 24.8}
\left[\begin{array}{cc}
T_{1} & T_{2} \\
0 & 0
\end{array}\right] 
\left[\begin{array}{cc}
T_{1}^{\dagger} & 0 \\
T_{2}^{\dagger} & 0
\end{array}\right]
\left[\begin{array}{cc}
T_{1} & T_{2} \\
0 & 0
\end{array}\right] = \left[\begin{array}{cc}
T_{1} & T_{2} \\
0 & 0
\end{array}\right] 
\left[\begin{array}{cc}
T_{1}^{\dagger}T_{1} & 0 \\
0 & T_{2}^{\dagger}T_{2}
\end{array}\right]= \left[\begin{array}{cc}
T_{1} & T_{2} \\
0 & 0
\end{array}\right].
\end{align}
\begin{align}\label{equ 25.8}
\left[\begin{array}{cc}
T_{1} & T_{2} \\
0 & 0
\end{array}\right] 
\left[\begin{array}{cc}
T_{1}^{\dagger} & 0 \\
T_{2}^{\dagger} & 0
\end{array}\right] = \left[\begin{array}{cc}
P_{R(T_{1})} + P_{R(T_{2})} & 0 \\
0 & 0
\end{array}\right] =  \left[\begin{array}{cc}
P_{R(T)} & 0 \\
0 & 0
\end{array}\right]
\end{align}
is symmetric. Moreover,
\begin{align}\label{equ 26.8}
\left[\begin{array}{cc}
T_{1}^{\dagger} & 0 \\
T_{2}^{\dagger} & 0
\end{array}\right]
\left[\begin{array}{cc}
T_{1} & T_{2} \\
0 & 0
\end{array}\right] = \left[\begin{array}{cc}
P_{N(T_{1})^{\perp}}\vert_{D(T_{1})} & 0 \\
0 & P_{N(T_{2})^{\perp}}\vert_{D(T_{2})}
\end{array}\right]
\end{align}
is also symmetric. Therefore, Theorem 5.7 in \cite{MR0451661} and the relations (\ref{equ 23.8}), (\ref{equ 24.8}), (\ref{equ 25.8}), (\ref{equ 26.8}) confirm that $T^{\dagger} =
\left[\begin{array}{cc}
T_{1}^{\dagger} & 0 \\
T_{2}^{\dagger} & 0
\end{array}\right]$.
\end{proof}
\begin{corollary}\label{cor 2.23.8}
\noindent Let $T =
\left[\begin{array}{cc}
T_{1} & T_{2} \\
0 & 0
\end{array}\right]$ be an operator on $H \bigoplus K$ with $T_{1} \in C(H)$ and $T_{2}\in C(K, H)$ both are densely defined closed range operators. If $R(T_{1}) \oplus^{\perp} R(T_{2}) \subset H$, then  $w(T) =
\left[\begin{array}{cc}
w(T_{1}) & w(T_{2}) \\
0 & 0
\end{array}\right]$.
\end{corollary}
\begin{proof}
$w(T) =(T^{\dagger})^{*}
 = \left[\begin{array}{cc}
T_{1}^{\dagger} & 0 \\
T_{2}^{\dagger} & 0
\end{array}\right]^{*} (\text{by Theorem \ref{thm 2.22.8}})$. Since $T_{1}^{\dagger}, T_{2}^{\dagger}$ both are bounded. Therefore, $w(T) =(T^{\dagger})^{*}
 = \left[\begin{array}{cc}
(T_{1}^{\dagger})^{*} & (T_{2}^{\dagger})^{*} \\
0 & 0
\end{array}\right] = \left[\begin{array}{cc}
w(T_{1}) & w(T_{2}) \\
0 & 0
\end{array}\right]$.
\end{proof}
\begin{corollary}\label{cor 2.24.8}
\noindent Let $T =
\left[\begin{array}{cc}
T_{1} & T_{2} \\
0 & 0
\end{array}\right]$ be an operator on $H \bigoplus K$ with $T_{1} \in C(H)$ and $T_{2}\in C(K, H)$ both are densely defined closed range operators. If $R(T_{1}) \oplus^{\perp} R(T_{2}) \subset H$, then  $w(\vert T \vert) =
\left[\begin{array}{cc}
w(\vert T_{1}\vert) & 0 \\
0 & w(\vert T_{2}\vert)
\end{array}\right]$.
\end{corollary}
\begin{proof}
 Theorem \ref{thm 2.22.8} and Corollary \ref{cor 2.9.8} guarantee $w(\vert T \vert) = \vert w(T) \vert$. By Corollary \ref{cor 2.23.8}, we get
\begin{align*}
{\vert w(T) \vert}^{2} = 
\left[\begin{array}{cc}
w(T_{1}) & w(T_{2}) \\
0 & 0
\end{array}\right]^{*} 
\left[\begin{array}{cc}
w(T_{1}) & w(T_{2}) \\
0 & 0
\end{array}\right] =
\left[\begin{array}{cc}
{\vert w(T_{1})\vert}^{2} & 0 \\
0 & {\vert w(T_{2}) \vert}^{2}
\end{array}\right].
\end{align*} 
Therefore,
\begin{align*}
w(\vert T \vert) = \left[\begin{array}{cc}
{\vert w(T_{1})\vert} & 0 \\
0 & {\vert w(T_{2}) \vert}
\end{array}\right] = \left[\begin{array}{cc}
{ w(\vert T_{1} \vert)} & 0 \\
0 & {w(\vert T_{2} \vert)}
\end{array}\right].
\end{align*}
\end{proof}
Now, we explore the Moore-Penrose inverse of the lower-triangular operator matrix $T = \left[\begin{array}{cc}
T_{1} & 0 \\
T_{3} & T_{4}
\end{array}\right]$, where $T_{1} \in C(H)$, $T_{3} \in C(H,K)$ and $T_{4} \in C(K)$, with all $T_{i} ~(i =1, 3, 4)$ being densely defined closed-range operators and $T$ satisfies  the following conditions:
\begin{enumerate}
\item $R(T_{1}^{*}) \oplus^{\perp} R(T_{3}^{*}) = H $,
\item $R(T_{3}) \oplus^{\perp} R(T_{4}) \subset K$.
\end{enumerate}
\begin{theorem}
When $T$ satisfies the above conditions, then $T^{\dagger}$ exists and 
$T^{\dagger} = \left[\begin{array}{cc}
T_{1}^{\dagger} & T_{3}^{\dagger} \\
0 & T_{4}^{\dagger}
\end{array}\right] = \left[\begin{array}{cc}
\overline{S^{\dagger}T_{1}^{*}} & \overline{S^{\dagger}T_{3}^{*}} \\
0 & T_{4}^{\dagger}
\end{array}\right]$, where $S = T_{1}^{*}T_{1} + T_{3}^{*}T_{3}$ with $D(S)$ is densely defined. Moreover,
$w(T) = \left[\begin{array}{cc}
w(T_{1}) & 0 \\
w(T_{3}) & w(T_{4})
\end{array}\right]$.
\end{theorem}
\begin{proof}
First, we claim that $S$ is closed with closed range. Let us consider $(x,y) \in \overline{G(S)}$, then there exists a sequence $\{x_{n}\} \in D(S)$ such that $x_{n} \to x$ and $(T_{1}^{*}T_{1} + T_{3}^{*}T_{3})x_{n} \to y$, as $n \to \infty$. So, $\{T_{1}^{*}T_{1}x_{n}\}$ and $\{T_{3}^{*}T_{3}x_{n}\}$ both are Cauchy sequences. We know that $T_{1}^{*}T_{1}$ and $T_{3}^{*}T_{3}$ are closed. So, $y = (T_{1}^{*}T_{1} + T_{3}^{*}T_{3})x$. Thus, $S$ is closed. For $v \in \overline{R(S)}$, there exists a sequence $\{u_{n}\}$ in $D(S)$ with $Su_{n} \to v$ as $n \to \infty$. Again, $\{T_{1}^{*}T_{1}u_{n}\}$ and $\{T_{3}^{*}T_{3}u_{n}\}$ both are Cauchy sequences. Now, we can write $u_{n} = u_{n}^{\prime} + u_{n}^{\prime{\prime}}$, where $u_{n}^{'} \in N(T_{1}^{*}T_{1})^{\perp} = R(T_{1}^{*})$ and $u_{n}^{\prime{\prime}} \in N(T_{1}^{*}T_{1}) = R(T_{1}^{*})^{\perp}, \text{ for all $n \in \mathbb{N}$}$. This implies that 
\begin{align*}
\|T_{1}^{*}T_{1}u_{n} - T_{1}^{*}T_{1}u_{m}\| \geq \gamma(T_{1}^{*}T_{1}) \|u_{n}^{\prime} - u_{m}^{\prime}\|
\end{align*}
\begin{align*}
\|T_{3}^{*}T_{3}u_{n} - T_{3}^{*}T_{3}u_{m}\| \geq \gamma(T_{3}^{*}T_{3}) \|u_{n}^{\prime{\prime}} - u_{m}^{\prime{\prime}}\|
\end{align*}

, where $\gamma({T_{1}^{*}T_{1}})$ and $\gamma({T_{3}^{*}T_{3}})$ are the reduced minimal modulus of $T_{1}^{*}T_{1}$ and $T_{3}^{*}T_{3}$, respectively. It confirms that $u_{n}^{\prime} \to u^{\prime}$ and $u_{n}^{\prime{\prime}} \to u^{\prime{\prime}}$ as $n \to \infty$, for some $u^{\prime} \in R(T_{1}^{*})$ and $u^{\prime{\prime}} \in R(T_{3}^{*})$. The closeness of $T_{1}^{*}T_{1}$ and $T_{3}^{*}T_{3}$ say that $T_{1}^{*}T_{1}u_{n} \to T_{1}^{*}T_{1}u$ and $T_{3}^{*}T_{3}u_{n} \to T_{3}^{*}T_{3}u$, as $n \to \infty$, where $u = u^{\prime} + u^{\prime{\prime}}$. Thus, $Su_{n} \to Su$, as $n \to \infty$ which implies $Su= v$. $R(S)$ is closed. It means that $S^{\dagger}$ exists and is bounded. Now, we will show that $S^{\dagger} = (T_{1}^{*}T_{1})^{\dagger} + (T_{3}^{*}T_{3})^{\dagger}$. 
\begin{align}\label{equ 27.8}
(T_{1}^{*}T_{1} + T_{3}^{*}T_{3})((T_{1}^{*}T_{1})^{\dagger} + (T_{3}^{*}T_{3})^{\dagger}) (T_{1}^{*}T_{1} + T_{3}^{*}T_{3}) = (T_{1}^{*}T_{1} + T_{3}^{*}T_{3}).
\end{align}
\begin{align}\label{equ 28.8}
((T_{1}^{*}T_{1})^{\dagger} + (T_{3}^{*}T_{3})^{\dagger}) (T_{1}^{*}T_{1} + T_{3}^{*}T_{3})((T_{1}^{*}T_{1})^{\dagger} + (T_{3}^{*}T_{3})^{\dagger}) = ((T_{1}^{*}T_{1})^{\dagger} + (T_{3}^{*}T_{3})^{\dagger}).
\end{align}
\begin{align}\label{equ 29.8}
((T_{1}^{*}T_{1})^{\dagger} + (T_{3}^{*}T_{3})^{\dagger}) (T_{1}^{*}T_{1} + T_{3}^{*}T_{3}) &= P_{N(T_{1})^{\perp}}\vert_{D(T_{1}^{*}T_{1})} + P_{N(T_{3})^{\perp}}\vert_{D(T_{3}^{*}T_{3})}\\
&= I_{D(S)} \subset I_{H}
\end{align}
is symmetric because $D(S)$ is dense.
\begin{align}\label{equ 30.8}
(T_{1}^{*}T_{1} + T_{3}^{*}T_{3})((T_{1}^{*}T_{1})^{\dagger} + (T_{3}^{*}T_{3})^{\dagger}) = P_{R(T_{1}^{*})} + P_{R(T_{3}^{*})} = I_{H}.
\end{align}
By Theorem 5.7 in \cite{MR0451661} and the above four relations (\ref{equ 27.8}), (\ref{equ 28.8}), (\ref{equ 29.8}), (\ref{equ 30.8}) guarantee that $S^{\dagger} = (T_{1}^{*}T_{1})^{\dagger} + (T_{3}^{*}T_{3})^{\dagger}$. Let us assume $\biggl( \begin{pmatrix} v_{1}\\v_{2}\end{pmatrix}, \begin{pmatrix}w_{1}\\w_{2} \end{pmatrix} \biggr) \in \overline{G(T)}$. There exists a sequence $\biggl\{\begin{pmatrix} v_{1n}\\ v_{2n} \end{pmatrix} \biggr\}$ such that $\begin{pmatrix} v_{1n}\\ v_{2n} \end{pmatrix} \to \begin{pmatrix} v_{1}\\ v_{2} \end{pmatrix}$ and $\begin{pmatrix} T_{1}v_{1n}\\ T_{3}v_{1n} + T_{4}v_{2n} \end{pmatrix} \to \begin{pmatrix} w_{1}\\ w_{2} \end{pmatrix}$ as $n \to \infty$. The given condition $R(T_{3}) \oplus^{\perp} R(T_{4})$ and the closeness of $T_{1}$, $T_{2}$, $T_{3}$ confirm that $\begin{pmatrix} T_{1}v_{1} \\ T_{3}v_{1} + T_{4}v_{2}\end{pmatrix} = \begin{pmatrix} w_{1}\\ w_{2} \end{pmatrix}$. Hence, $T$ is closed. It is easy to show that $R(T)$ is closed because $R(T_{1})$ and $R(T_{3}) \oplus^{\perp} R(T_{4})$ both are closed. Again,
\begin{align*}
S^{\dagger}T_{3}^{*} = {T_{3}}^{\dagger}\vert_{D(T_{3}^{*})} \subset T_{3}^{\dagger} \text{ and } T_{3}^{\dagger} = (T_{3}(S^{\dagger})^{*})^{*} = \overline{S^{\dagger}T_{3}^{*}}.
\end{align*}
Similarly, $\overline{S^{\dagger}T_{1}^{*}} = T_{1}^{\dagger}$. we will establish that $T^{\dagger} = \left[\begin{array}{cc}
T_{1}^{\dagger} & T_{3}^{\dagger} \\
0 & T_{4}^{\dagger}
\end{array}\right] $. 
\begin{align}\label{equ 32.8}
\left[\begin{array}{cc}
T_{1}^{\dagger} & T_{3}^{\dagger} \\
0 & T_{4}^{\dagger}
\end{array}\right]  \left[\begin{array}{cc}
T_{1} & 0 \\
T_{3} & T_{4}
\end{array}\right]  \left[\begin{array}{cc}
T_{1}^{\dagger} & T_{3}^{\dagger} \\
0 & T_{4}^{\dagger}
\end{array}\right] = \left[\begin{array}{cc}
T_{1}^{\dagger} & T_{3}^{\dagger} \\
0 & T_{4}^{\dagger}
\end{array}\right].
\end{align}
\begin{align}\label{equ 33.8}
\left[\begin{array}{cc}
T_{1} & 0 \\
T_{3} & T_{4}
\end{array}\right] \left[\begin{array}{cc}
T_{1}^{\dagger} & T_{3}^{\dagger} \\
0 & T_{4}^{\dagger}
\end{array}\right]  \left[\begin{array}{cc}
T_{1} & 0 \\
T_{3} & T_{4}
\end{array}\right] = \left[\begin{array}{cc}
T_{1} & 0 \\
T_{3} & T_{4}
\end{array}\right].
\end{align}
\begin{align}\label{35.8}
\left[\begin{array}{cc}
T_{1}^{\dagger} & T_{3}^{\dagger} \\
0 & T_{4}^{\dagger}
\end{array}\right]  \left[\begin{array}{cc}
T_{1} & 0 \\
T_{3} & T_{4}
\end{array}\right] &= \left[\begin{array}{cc}
P_{R(T_{1}^{*})}\vert_{D(T_{1})} + P_{R(T_{3}^{*})}\vert_{D(T_{3})} & 0 \\
0 & P_{R(T_{4}^{*})}\vert_{D(T_{4})}
\end{array}\right]  \\
&= \left[\begin{array}{cc}
I_{D(T_{1}) \cap D(T_{3})} & 0 \\
0 & P_{R(T_{4}^{*})}\vert_{D(T_{4})}
\end{array}\right]
\end{align}
is symmetric because $D(T_{1}) \cap D(T_{3}) \supset D(S)$ is dense.  Moreover,
\begin{align}\label{equ 36.8}
\left[\begin{array}{cc}
T_{1} & 0 \\
T_{3} & T_{4}
\end{array}\right] \left[\begin{array}{cc}
T_{1}^{\dagger} & T_{3}^{\dagger} \\
0 & T_{4}^{\dagger}
\end{array}\right] = 
\left[\begin{array}{cc}
P_{R(T_{1})} & 0 \\
0 & P_{R(T_{3})} + P_{R(T_{4})}
\end{array}\right] 
\end{align}
is symmetric. Therefore,
$T^{\dagger}= 
\left[\begin{array}{cc}
T_{1}^{\dagger} & T_{3}^{\dagger} \\
0 & T_{4}^{\dagger}
\end{array}\right] = \left[\begin{array}{cc}
\overline{S^{\dagger}T_{1}^{*}} & \overline{S^{\dagger}T_{3}^{*}}\\
0 & T_{4}^{\dagger}
\end{array}\right]$ and $w(T) = 
\left[\begin{array}{cc}
w(T_{1}) & 0 \\
w(T_{3}) & w(T_{4})
\end{array}\right]$.
\end{proof}
\begin{remark}
Let us consider $T = \left[\begin{array}{cc}
T_{1} & T_{2} \\
T_{3} & T_{4}
\end{array}\right]$, where $T_{1} \in C(H)$, $T_{2} \in C(K, H)$ $T_{3} \in C(H,K)$ and $T_{4} \in C(K)$, with $D(T_{1}) \cap D(T_{3})$ and $D(T_{2}) \cap D(T_{4})$ both are densely defined. Assume that all $R(T_{i}) ~(i= 1, 2, 3, 4)$ are closed. Moreover, $T$ satisfies  the following conditions:
\begin{enumerate}
\item $R(T_{1}^{*}) \oplus^{\perp} R(T_{3}^{*}) \subset H $,
\item $R(T_{2}^{*}) \oplus^{\perp} R(T_{4}^{*}) \subset K$,
\item $R(T_{1}) \oplus^{\perp} R(T_{2}) \subset H$,
\item $R(T_{3}) \oplus^{\perp} R(T_{4}) \subset K$.
\end{enumerate}
Then $T$ is closed and $T^{\dagger}$ is bounded. Moreover, $T^{\dagger}= 
\left[\begin{array}{cc}
T_{1}^{\dagger} & T_{3}^{\dagger} \\
T_{2}^{\dagger} & T_{4}^{\dagger}
\end{array}\right] $ and $w(T) = 
\left[\begin{array}{cc}
w(T_{1}) & w(T_{2}) \\
w(T_{3}) & w(T_{4})
\end{array}\right]$.

\end{remark}

\begin{center}
	\textbf{Acknowledgements}
\end{center}

\noindent The present work of the second author was partially supported by Science and Engineering Research Board (SERB), Department of Science and Technology, Government of India (Reference Number: MTR/2023/000471) under the scheme ``Mathematical Research Impact Centric Support (MATRICS)''. 

\section*{Declarations}
  The authors declare no conflicts of interest.

\end{document}